\documentclass[11pt]{article}

\usepackage[T1]{fontenc}
\usepackage{lmodern}
\usepackage[margin=1.1in]{geometry}
\usepackage{amsmath,amssymb,amsthm}
\usepackage{booktabs}
\usepackage{microtype}
\usepackage{tikz}
\usepackage[hidelinks]{hyperref}
\newtheorem{theorem}{Theorem}
\newtheorem{proposition}[theorem]{Proposition}
\newtheorem{lemma}[theorem]{Lemma}
\newtheorem{corollary}[theorem]{Corollary}
\theoremstyle{remark}
\newtheorem{remark}[theorem]{Remark}

\newcommand{\R}{\mathbb{R}}
\newcommand{\Z}{\mathbb{Z}}
\newcommand{\N}{\mathbb{N}}

\title{Lattice Trapezoids of Prime Gaps:\\
Minimal Interior Counts and Empirical Statistics}
\author{Rafael A. Gonzalez Castillo\\
\small Department of Mathematics, Lehman College (CUNY), Bronx, NY, USA}
\date{November 2026}

\begin{document}
\maketitle

\begin{abstract}
Let $p_n$ denote the $n$-th prime and $g_n = p_{n+1} - p_n$. Mapping consecutive primes
to $z_n = p_n + i\,p_{n+1} \in \Z[i]$ defines, for each $n \ge 1$, a convex lattice
quadrilateral $Q_n \subset \R^2$ whose vertices lie in $\Z^2$ with area $\Delta A_n = \tfrac12 g_n(g_n + g_{n+1})$.
We prove two-sided bounds relating $\sum_{n=1}^{N} \Delta A_n$ to the second moment
$\sum_{n=1}^{N} g_n^2$, express the boundary and interior lattice-point counts $B_n$ and
$I_n$ of $Q_n$ in terms of $(g_n, g_{n+1})$, and classify all minimal values of $I_n$. In
particular, for $n \ge 3$ one has $I_n \ge 2$, with equality if and only if
$(p_n, p_{n+1}, p_{n+2})$ is a prime triple $(p, p+2, p+6)$. We also prove a telescoping
identity for $\sum_{n=1}^{N} B_n$ implying
$\lim_{N\to\infty} \sum_{n=1}^{N} I_n / \sum_{n=1}^{N} \Delta A_n = 1$ at rate
$O(1/\ln N)$ using only Chebyshev's bound and Bertrand's postulate. Numerically, for all
primes below $3\times 10^{10}$ ($1.3\times 10^{9}$ primes) we compare the distribution of
$\gcd(g_n, g_{n+1})$ and the sign of the discrete convexity
$K_n = \tfrac12(g_n g_{n+2} - g_{n+1}^2)$ with a local-scale windowed empirical-marginal
null. The residual adjacency shifts decrease monotonically in magnitude across six disjoint
intervals from $10^5$ to $2\times 10^{10}$. A one-parameter $c/\ln x$ decay is adequate only
for the $K_n > 0$ shift and is rejected for the other three; two-parameter forms fit all
four, but the data cannot determine the asymptotic rate. The algebraic bounds,
minimal-count classifications, telescoping boundary sum, and direct discrete
interior-point formula are verified in Lean~4.
\end{abstract}



\section{Introduction}\label{sec:intro}

Let $\N = \{1, 2, 3, \dots\}$ be the set of positive integers, let $p_n$ denote the
$n$-th prime ($p_1 = 2$), and let $g_n = p_{n+1} - p_n$ be the $n$-th prime gap. The
second moment $S_N = \sum_{n=1}^{N} g_n^2$ has been studied since Cram\'er~\cite{cramer}.
Under the Riemann hypothesis, Selberg~\cite{selberg} proved that
$\sum_{p_n \le x} g_n^2 \ll x(\ln x)^3$. Unconditionally, Heath-Brown~\cite{hb} proved
$\sum_{p_n \le x} g_n^2 \ll x^{23/18+\varepsilon}$, Peck~\cite{peck} improved this to
$x^{5/4+\varepsilon}$, and Stadlmann~\cite{stadlmann} obtained the current best exponent
$x^{1.23+\varepsilon}$. The conjectured asymptotic $\sum_{p_n \le x} g_n^2 \sim 2x\ln x$
and its numerical behavior are discussed in~\cite{wolf,cohen}.

For each $n \in \N$, consider the Gaussian integer $z_n = p_n + i\,p_{n+1} \in \Z[i]$.
The region in $\R^2$ bounded above by the line segment $z_n z_{n+1}$, below by the
diagonal $y = x$, and on the left and right by the vertical segments $x = p_n$ and
$x = p_{n+1}$ is a convex lattice quadrilateral $Q_n \subset \R^2$ whose vertices are
$(p_n, p_n)$, $(p_{n+1}, p_{n+1})$, $(p_{n+1}, p_{n+2})$, and $(p_n, p_{n+1})$. Its area
is $\Delta A_n = \tfrac12 g_n(g_n + g_{n+1})$. As an identity for real sequences, this
relates $\sum_{n=1}^{N} \Delta A_n$ to $S_N$ up to a factor of~$2$
(Section~\ref{sec:construction}).

Because the vertices of $Q_n$ lie in $\Z^2$, Pick's theorem expresses $\Delta A_n$ in
terms of the boundary lattice-point count $B_n$ and the strict interior lattice-point
count $I_n$ of $Q_n$. In Section~\ref{sec:pick}, we show that
$B_n = 2g_n + g_{n+1} + \gcd(g_n, g_{n+1})$, classify all minimal values of $I_n$, and
prove unconditionally that $\sum_{n=1}^{N} I_n / \sum_{n=1}^{N} \Delta A_n = 1 - O(1/\ln N)$.
Section~\ref{sec:numerics} reports the lattice invariants for primes up to $2\times 10^8$
and, for primes below $3\times 10^{10}$, separates single-gap marginal effects from
adjacent-gap shifts for both $\gcd(g_n, g_{n+1})$ and the discrete convexity
$K_n = \tfrac12(g_n g_{n+2} - g_{n+1}^2)$ using a windowed local-scale null, with
block-bootstrap uncertainties across six disjoint intervals. Section~\ref{sec:lean}
describes the Lean~4 verification of the finite and algebraic results in
Sections~\ref{sec:construction}--\ref{sec:pick}, including a direct lattice-point
verification of the interior-point formula for $Q_n$.

The geometric formulas encode adjacent prime gaps rather than imposing new
constraints on their occurrence. The contributions here are the explicit
lattice-point classifications, a telescoping boundary identity, a finite Lean~4
formalization, and empirical comparisons with a local marginal null model.
No new existence result for prime tuples or asymptotic law for the residual
adjacency shifts is asserted.

\section{The trapezoid construction and area bounds}\label{sec:construction}

Let $z_n = p_n + i\,p_{n+1}$ for $n \in \N$. Then $\mathrm{Re}(z_{n+1}) = \mathrm{Im}(z_n)
= p_{n+1}$ and $\Delta z_n = z_{n+1} - z_n = g_n + i\,g_{n+1}$. Since $g_n / p_n \to 0$,
we have $\arg(z_n) = \arctan(1 + g_n/p_n) \to \pi/4$, so the polyline $(z_n)$ approaches
the diagonal $y = x$ from above.

Dropping perpendiculars from $z_n$ and $z_{n+1}$ to the real axis gives a right trapezoid
of width $g_n$ and heights $p_{n+1}, p_{n+2}$, with area
\begin{equation}\label{eq:A}
A_n = \tfrac12\,(p_{n+1} + p_{n+2})(p_{n+1} - p_n).
\end{equation}
The corresponding trapezoid under $y = x$ over $[p_n, p_{n+1}]$ has area
\begin{equation}\label{eq:E}
E_n = \tfrac12\,(p_{n+1} + p_n)(p_{n+1} - p_n) = \tfrac12\,(p_{n+1}^2 - p_n^2),
\end{equation}
and we define $\Delta A_n = A_n - E_n$.

\begin{lemma}\label{lem:identity}
For all $n \in \N$,
\[
\Delta A_n = \tfrac12\,g_n\,(p_{n+2} - p_n) = \tfrac12\,g_n\,(g_n + g_{n+1}).
\]
Consequently, for every $N \in \N$,
\begin{equation}\label{eq:sum}
\sum_{n=1}^{N} \Delta A_n
  = \tfrac12 \sum_{n=1}^{N} g_n^2 + \tfrac12 \sum_{n=1}^{N} g_n g_{n+1}.
\end{equation}
\end{lemma}
\begin{proof}
Subtracting \eqref{eq:E} from \eqref{eq:A} yields
\[
\Delta A_n = \tfrac12 g_n\bigl[(p_{n+1}+p_{n+2})-(p_{n+1}+p_n)\bigr]
= \tfrac12 g_n(p_{n+2}-p_n) = \tfrac12 g_n(g_n + g_{n+1}).
\]
Summing over $1 \le n \le N$ gives \eqref{eq:sum}. This completes the proof.
\end{proof}

Note that Lemma~\ref{lem:identity} holds for any real sequence $(q_n)$, with
$\Delta A_n > 0$ whenever $(q_n)$ is strictly increasing. For primes, parity yields
integrality for $n \ge 2$.

\begin{proposition}\label{prop:basic}
For all $n \in \N$, $\Delta A_n > 0$. For $n \ge 2$, the gaps $g_n = 2a$ and
$g_{n+1} = 2b$ ($a,b \in \N$) are even, so $\Delta A_n = 2a(a+b)$ is an even integer with
$\Delta A_n \ge 4$ (with equality at $n=2$). At $n=1$, $\Delta A_1 = \tfrac32$.
\end{proposition}
\begin{proof}
For $n \ge 2$, $p_n$ is odd, so consecutive differences are even and
$\Delta A_n = \tfrac12(2a)(2a+2b) = 2a(a+b) \ge 4$. Direct evaluation gives
$\Delta A_1 = \tfrac12(1)(1+2) = \tfrac32$ and $\Delta A_2 = \tfrac12(2)(2+2) = 4$. This
completes the proof.
\end{proof}

Write $S_N = \sum_{n=1}^{N} g_n^2$.

\begin{theorem}[Sandwich and Cauchy--Schwarz bounds]\label{thm:sandwich}
For every $N \in \N$,
\begin{align}
\tfrac12\,S_N \;\le\; \sum_{n=1}^{N} \Delta A_n \;&\le\; S_{N+1}, \label{eq:sandwich}\\
\sum_{n=1}^{N} \Delta A_n \;&\ge\; \frac{(p_{N+1} - 2)^2}{2N}. \label{eq:lower}
\end{align}
\end{theorem}
\begin{proof}
The left inequality in \eqref{eq:sandwich} follows from \eqref{eq:sum} since
$g_n g_{n+1} > 0$. For the right inequality, $g_n g_{n+1} \le \tfrac12(g_n^2 + g_{n+1}^2)$
gives
\[
\sum_{n=1}^{N} g_n g_{n+1} \le \tfrac12 S_N + \tfrac12(S_{N+1} - g_1^2),
\]
whence $\sum_{n=1}^{N} \Delta A_n \le \tfrac34 S_N + \tfrac14 S_{N+1} \le S_{N+1}$.
For \eqref{eq:lower}, $\sum_{n=1}^{N} g_n = p_{N+1} - 2$ telescopes, and Cauchy--Schwarz
gives $(p_{N+1} - 2)^2 \le N S_N$; combining with the left inequality in
\eqref{eq:sandwich} yields \eqref{eq:lower}. This completes the proof.
\end{proof}

\begin{corollary}\label{cor:pnt}
By the prime number theorem ($p_{N+1} \sim N \ln N$), \eqref{eq:lower} implies
\[
\sum_{n=1}^{N} \Delta A_n \ge \bigl(\tfrac12 + o(1)\bigr)\,N(\ln N)^2.
\]
The comparison $\tfrac12 S_N \le \sum_{n=1}^{N}\Delta A_n \le S_{N+1}$
transfers estimates between the two sums, with the indicated index shift.
\end{corollary}

\section{Lattice geometry and Pick's theorem}\label{sec:pick}

Because $z_n = p_n + i\,p_{n+1} \in \Z[i]$, the region $Q_n \subset \R^2$ is a convex
lattice quadrilateral with vertices
\[
(p_n, p_n),\quad (p_{n+1}, p_{n+1}),\quad (p_{n+1}, p_{n+2}),\quad (p_n, p_{n+1}).
\]
Applying Pick's theorem to $Q_n$ expresses $\Delta A_n = \mathrm{Area}(Q_n)$ in terms of
its boundary lattice-point count $B_n$ and strict interior lattice-point count $I_n$.

\begin{proposition}[Pick decomposition]\label{prop:pick}
For every $n \in \N$, the boundary and interior lattice-point counts $B_n$ and $I_n$ of
$Q_n$ are given by
\begin{align}
B_n &= 2g_n + g_{n+1} + \gcd(g_n, g_{n+1}), \label{eq:Bn}\\
I_n &= \tfrac12\,g_n(g_n + g_{n+1}) - g_n - \tfrac12\,g_{n+1}
       - \tfrac12\,\gcd(g_n, g_{n+1}) + 1. \label{eq:In}
\end{align}
\end{proposition}
\begin{proof}
A lattice segment with displacement vector $(u,v) \in \Z^2$ contains $\gcd(|u|,|v|)$
primitive lattice steps. The four edges of $Q_n$ have displacement vectors $(g_n,g_n)$,
$(0,g_{n+1})$, $(-g_n,-g_{n+1})$, and $(0,-g_n)$, contributing $g_n$, $g_{n+1}$,
$\gcd(g_n,g_{n+1})$, and $g_n$ boundary steps, respectively. Summing gives \eqref{eq:Bn},
and Pick's theorem ($\Delta A_n = I_n + \tfrac12 B_n - 1$) together with
Lemma~\ref{lem:identity} gives \eqref{eq:In}. This completes the proof.
\end{proof}

\begin{remark}
Because $\Delta A_n$, $B_n$, and $I_n$ are functions of the adjacent gap pair
$(g_n, g_{n+1})$, any classification of $Q_n$ by its lattice invariants is equivalent to
a classification of the pairs $(g_n, g_{n+1})$. Pick's theorem therefore provides a
lattice-point reformulation of adjacent prime gaps rather than an independent arithmetic
constraint.
\end{remark}

\begin{theorem}[Classification of minimal interior counts]\label{thm:min_I}
For all $n \in \N$:
\begin{enumerate}
\item[(a)] $I_n = 0$ if and only if $n = 1$ (corresponding to $p_1 = 2$).
\item[(b)] $I_n = 1$ if and only if $n = 2$ (corresponding to $(p_2,p_3,p_4) = (3,5,7)$).
\item[(c)] If $n \ge 2$ and $g_n = 2$, then $I_n = g_{n+1}/2$. In particular, for
  $n \ge 3$ one has $I_n \ge 2$, with $I_n = 2$ if and only if $(g_n, g_{n+1}) = (2,4)$,
  that is, $(p_n, p_{n+1}, p_{n+2})$ is a prime triple of the form $(p, p+2, p+6)$.
\item[(d)] If $g_n \ge 4$, then $I_n \ge 7$, with equality if and only if
  $(g_n, g_{n+1}) = (4,2)$.
\end{enumerate}
\end{theorem}
\begin{proof}
Direct evaluation gives $I_1 = 0$ and $I_2 = 1$. Now assume $n \ge 2$, so $g_n$ and
$g_{n+1}$ are even and at least~$2$. Rewriting \eqref{eq:In} gives
\begin{equation}\label{eq:2In}
2I_n = (g_n - 2)(g_n + g_{n+1}) + 2 + \bigl(g_{n+1} - \gcd(g_n, g_{n+1})\bigr).
\end{equation}
Since $\gcd(g_n, g_{n+1}) \le g_{n+1}$, the final term is nonnegative.

If $g_n = 2$, then $\gcd(2, g_{n+1}) = 2$, and \eqref{eq:2In} reduces to
$2I_n = g_{n+1}$, so $I_n = g_{n+1}/2$. For $n \ge 3$ we have $p_n \ge 5$; since
$p_n, p_n+2, p_n+4$ form a complete residue system modulo~$3$, one of them is divisible
by~$3$ and hence composite. Thus $(g_n, g_{n+1}) \ne (2,2)$, so $g_{n+1} \ge 4$ and
$I_n \ge 2$, with equality if and only if $g_{n+1} = 4$.

If $g_n \ge 4$, then $(g_n - 2)(g_n + g_{n+1}) + 2 \ge 2(4+2) + 2 = 14$, so $I_n \ge 7$.
Equality holds if and only if $g_n = 4$, $g_{n+1} = 2$, and $\gcd(4,2) = 2$. Combining
these cases establishes (a)--(d). This completes the proof.
\end{proof}

\begin{figure}[htbp]
\centering
\begin{tikzpicture}[scale=0.68]
  \draw[step=1, gray!25, very thin] (4.5,4.5) grid (7.5,11.5);
  \draw[dashed, gray!70, thick] (4.5,4.5) -- (7.5,7.5)
    node[below right, font=\scriptsize, black!70] {$y=x$};
  \fill[blue!10] (5,5) -- (7,7) -- (7,11) -- (5,7) -- cycle;
  \draw[thick, blue!70!black] (5,5) -- (7,7) -- (7,11) -- (5,7) -- cycle;
  \foreach \x/\y in {5/5, 6/6, 7/7, 7/8, 7/9, 7/10, 7/11, 6/9, 5/7, 5/6} {
    \filldraw[fill=white, draw=blue!70!black, thick] (\x,\y) circle (2.8pt);
  }
  \foreach \x/\y in {6/7, 6/8} {
    \fill[red!80!black] (\x,\y) circle (3.2pt);
  }
  \node[below left, font=\scriptsize] at (5,5) {$(5,5)$};
  \node[below right, font=\scriptsize] at (7,7) {$(7,7)$};
  \node[above right, font=\scriptsize] at (7,11) {$z_4=(7,11)$};
  \node[above left, font=\scriptsize] at (5,7) {$z_3=(5,7)$};
  \node[below, font=\small] at (6,4.1) {$Q_3$: $(g_3,g_4)=(2,4)$};
  \node[below, font=\footnotesize] at (6,3.4) {$\Delta A_3 = 6,\; B_3 = 10,\; I_3 = 2$};
\end{tikzpicture}
\hspace{1.4cm}
\begin{tikzpicture}[scale=0.68]
  \draw[step=1, gray!25, very thin] (6.5,6.5) grid (11.5,13.5);
  \draw[dashed, gray!70, thick] (6.5,6.5) -- (11.5,11.5)
    node[below right, font=\scriptsize, black!70] {$y=x$};
  \fill[blue!10] (7,7) -- (11,11) -- (11,13) -- (7,11) -- cycle;
  \draw[thick, blue!70!black] (7,7) -- (11,11) -- (11,13) -- (7,11) -- cycle;
  \foreach \x/\y in {7/7, 8/8, 9/9, 10/10, 11/11, 11/12, 11/13, 9/12, 7/11, 7/10, 7/9, 7/8} {
    \filldraw[fill=white, draw=blue!70!black, thick] (\x,\y) circle (2.8pt);
  }
  \foreach \x/\y in {8/9, 8/10, 8/11, 9/10, 9/11, 10/11, 10/12} {
    \fill[red!80!black] (\x,\y) circle (3.2pt);
  }
  \node[below left, font=\scriptsize] at (7,7) {$(7,7)$};
  \node[below right, font=\scriptsize] at (11,11) {$(11,11)$};
  \node[above right, font=\scriptsize] at (11,13) {$z_5=(11,13)$};
  \node[above left, font=\scriptsize] at (7,11) {$z_4=(7,11)$};
  \node[below, font=\small] at (9,6.1) {$Q_4$: $(g_4,g_5)=(4,2)$};
  \node[below, font=\footnotesize] at (9,5.4) {$\Delta A_4 = 12,\; B_4 = 12,\; I_4 = 7$};
\end{tikzpicture}
\caption{The lattice trapezoids $Q_3$ (left) and $Q_4$ (right) bounded by the prime
segment $z_n z_{n+1}$ and the diagonal $y=x$. Hollow circles mark boundary lattice
points ($B_n$) and solid circles mark strict interior lattice points ($I_n$).}\label{fig:trapezoids}
\end{figure}
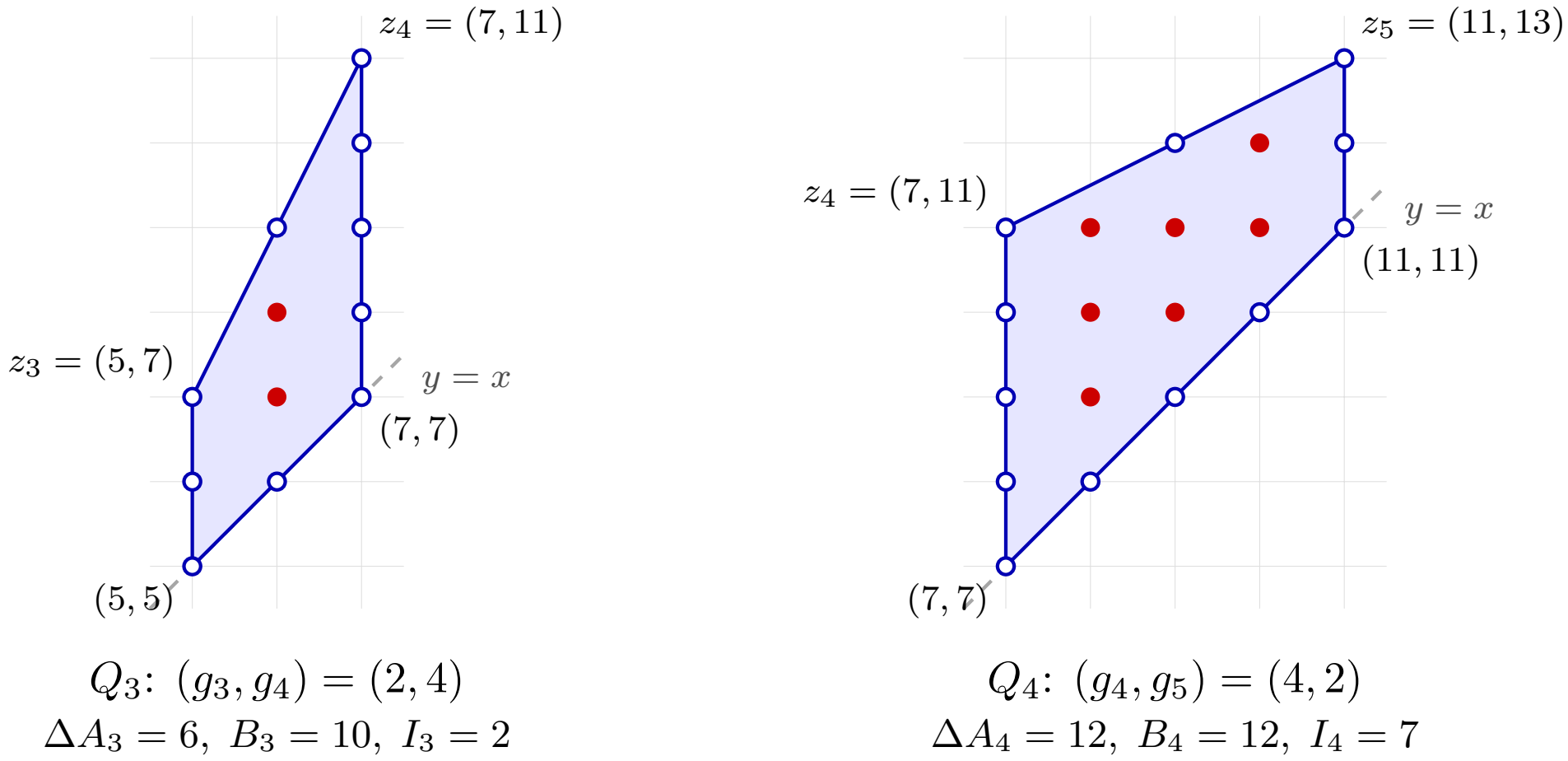

Figure~\ref{fig:trapezoids} illustrates the sharp minima of Theorem~\ref{thm:min_I}(c)
and (d) at $Q_3$ and $Q_4$.

\begin{remark}
Any prime triple $(p, p+2, p+6)$ with $p \ge 5$ consists of consecutive primes:
since $p$ and $p+2$ are prime and exceed $3$, one has $p \equiv 2 \pmod 3$,
$p+4 \equiv 0 \pmod 3$, and there is no intervening prime. By Theorem~\ref{thm:min_I}(c), the special case of the
Hardy--Littlewood prime $k$-tuples conjecture asserting that there are infinitely many
triples $(p, p+2, p+6)$ is equivalent to the statement that $I_n = 2$ occurs infinitely
often.
\end{remark}

\begin{theorem}[Telescoping boundary and interior density]\label{thm:boundary}
For every $N \in \N$,
\begin{equation}\label{eq:Bsum}
\sum_{n=1}^{N} B_n = 3p_{N+1} + g_{N+1} - 7 + \sum_{n=1}^{N} \gcd(g_n, g_{n+1})
\;\le\; 4p_{N+1} + g_{N+1} - 9.
\end{equation}
Consequently,
\[
\lim_{N \to \infty} \frac{\sum_{n=1}^{N} I_n}{\sum_{n=1}^{N} \Delta A_n} = 1,
\]
with $1 - \bigl(\sum_{n=1}^{N} I_n\bigr)/\bigl(\sum_{n=1}^{N} \Delta A_n\bigr)
= O\bigl(1/\ln N\bigr)$.
\end{theorem}
\begin{proof}
Summing \eqref{eq:Bn} and telescoping $\sum_{n=1}^{N} g_n = p_{N+1} - 2$ and
$\sum_{n=1}^{N} g_{n+1} = p_{N+2} - 3 = p_{N+1} + g_{N+1} - 3$ yields the identity in
\eqref{eq:Bsum}. The upper bound follows from $\gcd(g_n, g_{n+1}) \le g_n$. Combining
\eqref{eq:Bsum} with the lower bound \eqref{eq:lower} gives
\[
\frac{\tfrac12 \sum_{n=1}^{N} B_n}{\sum_{n=1}^{N} \Delta A_n}
  = O\!\left(\frac{N(p_{N+1} + g_{N+1})}{p_{N+1}^2}\right)
  = O\!\left(\frac{N}{p_{N+1}}\right),
\]
since $g_{N+1} \le p_{N+1}$ by Bertrand's postulate. By Chebyshev's elementary bound
$p_{N+1} \gg N \ln N$, this ratio is $O(1/\ln N)$. Dividing
$\sum_{n=1}^{N} I_n = \sum_{n=1}^{N} \Delta A_n - \tfrac12\sum_{n=1}^{N} B_n + N$ by
$\sum_{n=1}^{N} \Delta A_n$ yields the limit and the $O(1/\ln N)$ rate. (The same
argument shows that $\sum I_n / \sum \Delta A_n \to 1$ for any strictly increasing
integer sequence with $p_N / N \to \infty$ and $g_{N+1} = o(p_{N+1}^2 / N)$.) This
completes the proof.
\end{proof}

\section{Numerical experiments}\label{sec:numerics}

\paragraph{Lattice invariants (primes up to $2\times 10^8$).} We sieved all primes up to
$2\times 10^8$ ($11{,}078{,}937$ primes) and computed $\Delta A_n$, $B_n$, $I_n$, and
$\gcd(g_n, g_{n+1})$ in exact 64-bit integer arithmetic. Table~\ref{tab:small} records
the first ten prime trapezoids, illustrating Theorem~\ref{thm:min_I}: $I_1 = 0$,
$I_2 = 1$, and $I_n = 2$ occurring at the prime triples $(5,7,11)$, $(11,13,17)$, and
$(17,19,23)$.

\begin{table}[htbp]
\centering
\begin{tabular}{rrrrrrr}
\toprule
$n$ & $p_n$ & $(g_n, g_{n+1})$ & $\Delta A_n$ & $B_n$ & $\gcd(g_n, g_{n+1})$ & $I_n$ \\
\midrule
1 & 2 & $(1,2)$ & $1.5$ & 5 & 1 & 0 \\
2 & 3 & $(2,2)$ & $4.0$ & 8 & 2 & 1 \\
3 & 5 & $(2,4)$ & $6.0$ & 10 & 2 & 2 \\
4 & 7 & $(4,2)$ & $12.0$ & 12 & 2 & 7 \\
5 & 11 & $(2,4)$ & $6.0$ & 10 & 2 & 2 \\
6 & 13 & $(4,2)$ & $12.0$ & 12 & 2 & 7 \\
7 & 17 & $(2,4)$ & $6.0$ & 10 & 2 & 2 \\
8 & 19 & $(4,6)$ & $20.0$ & 16 & 2 & 13 \\
9 & 23 & $(6,2)$ & $24.0$ & 16 & 2 & 17 \\
10 & 29 & $(2,6)$ & $8.0$ & 12 & 2 & 3 \\
\bottomrule
\end{tabular}
\caption{Lattice invariants of the first ten prime trapezoids $Q_n$.}\label{tab:small}
\end{table}

Table~\ref{tab:density} tracks the cumulative interior-to-area ratio
$\sum_{n=1}^{N} I_n / \sum_{n=1}^{N} \Delta A_n$ and the boundary contribution
$\sum_{n=1}^{N} (B_n/2) / \sum_{n=1}^{N} \Delta A_n$ across logarithmic checkpoints up
to $2\times 10^8$. The interior ratio rises monotonically from $0.8463$ at $x = 10^4$ to
$0.9366$ at $x = 2\times 10^8$, while the boundary ratio decays from $0.1659$ to
$0.0656$, consistent with the $O(1/\ln N)$ rate of Theorem~\ref{thm:boundary}.

\begin{table}[htbp]
\centering
\begin{tabular}{rrrr}
\toprule
$x$ & $N$ & $\sum I_n / \sum \Delta A_n$ & $\sum (B_n/2) / \sum \Delta A_n$ \\
\midrule
$10^{4}$ & 1{,}229 & 0.846279 & 0.165865 \\
$10^{5}$ & 9{,}592 & 0.882977 & 0.124236 \\
$10^{6}$ & 78{,}498 & 0.906900 & 0.097798 \\
$10^{7}$ & 664{,}579 & 0.922651 & 0.080656 \\
$10^{8}$ & 5{,}761{,}455 & 0.933837 & 0.068618 \\
$2\times 10^8$ & 11{,}078{,}935 & 0.936633 & 0.065628 \\
\bottomrule
\end{tabular}
\caption{Convergence of the interior-to-area ratio for prime trapezoids up to
$2\times 10^8$ ($N$ denotes the number of trapezoids with $p_n \le x$, capped at
$\pi(2\times 10^8)-2$ at the final endpoint).}\label{tab:density}
\end{table}

\paragraph{Gap statistics (primes below $3\times 10^{10}$).} The remaining statistics
concern only the gaps, so they can be computed in a single streaming pass that never
stores the primes. A segmented sieve generates the primes in chunks, and for every
anchor $p_n \ge 3$ we use the triple $(g_n, g_{n+1}, g_{n+2})$ to accumulate the
histogram of $g_n$, the counts of $\gcd(g_n, g_{n+1})$, and the counts of the sign of
$K_n$, per multiplicative window and per bootstrap block (below). The sieve found
$1{,}300{,}005{,}926$ primes below $3\times 10^{10}$, matching the known value of
$\pi(3\times 10^{10})$, with largest gap $456$; this yields $N = 1{,}300{,}005{,}922$
anchors. (At $2\times 10^{8}$ the same code reproduces $\pi(2\times 10^8) = 11{,}078{,}937$
and the $11{,}078{,}933$ anchors of the earlier computation.)

\paragraph{Windowed marginal null.} Because prime gaps grow with $x$ (mean $\sim \ln x$),
pooling gaps across the whole range mixes scales. To control for local scale dependence,
we construct a \emph{windowed empirical-marginal null} by partitioning $[3, x_{\max}]$ into
an initial bin $[3, 10^3)$ (containing only $167$ primes) followed by multiplicative windows
$[t, 1.1t)$, computing exact independent-draw probabilities from the empirical single-gap
histogram inside each window, and aggregating weighted by window size. On
$[3, 2\times 10^8]$, pooled and windowed nulls differ by about $0.02$ percentage points
for convexity and by about $0.01$ percentage points for $\gcd(g_n, g_{n+1})$. Repeating the full-range computation with wider ratios $1.2$ and $1.5$ yields
similar residual shifts. Relative to the $1.1$-window results in
Table~\ref{tab:scale}, the displayed shifts change by at most approximately $0.013$
percentage points (for $\gcd=2$ in the smallest interval); in intervals beginning
at $10^7$ or above, the differences are at most $0.002$ percentage points at the
reported precision. All four shifts keep their signs and decrease in absolute
magnitude across the six intervals for each of the three ratios. The
one-parameter fits also retain the same qualitative goodness-of-fit conclusions
under the two wider window ratios. These finite-range sensitivity checks do not
identify an asymptotic decay law.

\paragraph{Distribution of $\gcd(g_n, g_{n+1})$.} For even $m \ge 2$, the divisibility
condition $m \mid \gcd(g_n, g_{n+1})$ holds if and only if
$p_n \equiv p_{n+1} \equiv p_{n+2} \pmod m$. The distribution of $\gcd(g_n, g_{n+1})$
therefore reformulates consecutive-prime residue statistics of the type studied by
Lemke~Oliver and Soundararajan~\cite{los}. Table~\ref{tab:gcd} compares the empirical
frequencies for $p_n \ge 3$, $p_n < 3\times 10^{10}$ (mean $\gcd = 4.21977$) against the
windowed empirical-marginal null and the \emph{limiting reduced-residue model}, in which
vanishing residue biases yield $\mathbb{P}(m \mid \gcd(g_n, g_{n+1})) = 1/\varphi(m)^2$
and, by M\"obius inversion,
\[
d_{2k} = \sum_{j=1}^{\infty} \frac{\mu(j)}{\varphi(2kj)^2}
\]
(evaluated numerically by truncating at $j = 5000$).

\begin{table}[htbp]
\centering
\begin{tabular}{rrrrr}
\toprule
$\gcd(g_n, g_{n+1})$ & Count & Empirical & Windowed null & Reduced-residue limit \\
\midrule
2 & 775{,}242{,}957 & $59.634\%$ & $58.312\%$ & $49.51\%$ \\
4 & 174{,}176{,}234 & $13.398\%$ & $13.367\%$ & $12.38\%$ \\
6 & 197{,}582{,}665 & $15.199\%$ & $15.336\%$ & $14.67\%$ \\
8 & 34{,}517{,}148 & $2.655\%$ & $2.741\%$ & $3.09\%$ \\
10 & 33{,}691{,}162 & $2.592\%$ & $2.748\%$ & $3.17\%$ \\
12 & 36{,}588{,}292 & $2.814\%$ & $3.016\%$ & $3.67\%$ \\
\bottomrule
\end{tabular}
\caption{Distribution of $\gcd(g_n, g_{n+1})$ for $3 \le p_n < 3\times 10^{10}$
($1{,}300{,}005{,}922$ anchors), compared with the windowed empirical-marginal null and
the limiting reduced-residue prediction $d_{2k}$.}\label{tab:gcd}
\end{table}

Table~\ref{tab:gcd} separates single-gap residue bias from adjacent-gap correlation.
Both baselines predict the orderings $\gcd = 6 > \gcd = 4$ and
$\gcd = 12 > \gcd = 8, 10$; the windowed null ($15.34\%$ vs.\ $13.37\%$) reflects the
excess of gap~$6$ among single gaps (the jumping champion in this range~\cite{orw}).
Of the $10.12$ percentage-point gap between the empirical frequency at $\gcd = 2$
($59.634\%$) and the limiting reduced-residue model ($49.51\%$), the windowed marginal
null accounts for $8.80$ points (about $87\%$), while $1.32$ points (about $13\%$)
reflects adjacency. Similarly, $\mathbb{P}(3 \mid \gcd) = 0.19798$ (vs.\ $0.20429$
windowed null and $0.25000$ limit) and $\mathbb{P}(4 \mid \gcd) = 0.20144$ (vs.\
$0.20727$ windowed null and $0.25000$ limit): three consecutive primes sharing a residue
class occur less often than local independence predicts, consistent in sign with the
Lemke~Oliver--Soundararajan bias against repeated residues~\cite{los}. Both the
single-gap deficit and the adjacency shift are smaller than at $2\times 10^{8}$,
where $\mathbb{P}(3 \mid \gcd) = 0.18600$, in line with slowly decaying biases.

\paragraph{Discrete convexity.} For $n \ge 2$, the signed area of the lattice triangle
with vertices $z_n, z_{n+1}, z_{n+2}$ is $K_n = \tfrac12(g_n g_{n+2} - g_{n+1}^2)$.
Because slopes increase along $(z_n)$ precisely when $g_{n+1}/g_n < g_{n+2}/g_{n+1}$,
$K_n > 0$ represents local convexity of the polyline as a graph. Across
$N = 1{,}300{,}005{,}922$ triples below $3\times 10^{10}$ we observe $K_n > 0$ in
$622{,}969{,}431$ cases ($47.921\%$ vs.\ $48.119\%$ windowed null, shift
$-0.198\text{ pp}$), $K_n = 0$ in $7{,}016{,}616$ cases ($0.540\%$ vs.\ $0.681\%$, shift
$-0.141\text{ pp}$), and $K_n < 0$ in $670{,}019{,}875$ cases ($51.540\%$ vs.\
$51.200\%$, shift $+0.340\text{ pp}$). Most of the asymmetry between concavity and
convexity is explained by the skewness of the single-gap marginal alone.

\paragraph{Disjoint-window scale dependence.} To test how the adjacency shifts depend on
scale without cumulative averaging, Table~\ref{tab:scale} evaluates the shifts (empirical
minus the local $1.1\times$-windowed marginal null) on six disjoint intervals. Because
consecutive gap pairs and triples overlap, standard errors are computed by a
non-overlapping block bootstrap: each interval is split into $200$ contiguous blocks of
equal width in $x$ (about $82$ primes per block in the smallest interval and about
$2.1\times 10^{6}$ in the largest), and $1{,}000$ resamples of the blocks are drawn, with
the windowed null recomputed from the resampled gap histograms in each resample. These
standard errors quantify sampling error only.

\begin{table}[htbp]
\centering
\small
\begin{tabular}{lrrrrr}
\toprule
Interval & $N$ & $\gcd=2$ shift & $K_n > 0$ shift & $K_n = 0$ shift & $K_n < 0$ shift \\
\midrule
$[10^5, 3\times 10^5]$ & 16{,}405 & $+2.902 \pm 0.243$ & $-0.480 \pm 0.215$ & $-0.360 \pm 0.097$ & $+0.840 \pm 0.209$ \\
$[10^6, 3\times 10^6]$ & 138{,}318 & $+2.241 \pm 0.088$ & $-0.314 \pm 0.060$ & $-0.333 \pm 0.030$ & $+0.647 \pm 0.065$ \\
$[10^7, 3\times 10^7]$ & 1{,}193{,}280 & $+2.006 \pm 0.030$ & $-0.273 \pm 0.022$ & $-0.242 \pm 0.008$ & $+0.514 \pm 0.022$ \\
$[10^8, 2\times 10^8]$ & 5{,}317{,}482 & $+1.708 \pm 0.016$ & $-0.239 \pm 0.011$ & $-0.204 \pm 0.004$ & $+0.443 \pm 0.011$ \\
$[10^9, 3\times 10^9]$ & 93{,}602{,}003 & $+1.446 \pm 0.004$ & $-0.213 \pm 0.002$ & $-0.161 \pm 0.001$ & $+0.374 \pm 0.002$ \\
$[10^{10}, 2\times 10^{10}]$ & 427{,}154{,}205 & $+1.293 \pm 0.002$ & $-0.196 \pm 0.001$ & $-0.137 \pm 0.000$ & $+0.332 \pm 0.001$ \\
\bottomrule
\end{tabular}
\caption{Residual shifts (empirical minus local windowed marginal null, in percentage
points) with block-bootstrap standard errors (null recomputed within each resample)
across six disjoint intervals. Standard errors below $0.0005$ are shown as $0.000$.}\label{tab:scale}
\end{table}

All four shifts shrink in magnitude at every step across the six intervals. The decrease
between the last two intervals exceeds seven standard errors in every column (about
$34\sigma$, $8\sigma$, $20\sigma$, and $19\sigma$ for $\gcd=2$, $K_n>0$, $K_n=0$, and
$K_n<0$, computed from the rounded standard errors in the table); earlier steps are less
significant, especially for $K_n > 0$ ($0.6\sigma$ to $2.3\sigma$ per step before the
last).

We fitted each shift against $\ln x$, taken as the interval mean of $\ln p_n$,
by inverse-variance weighted least squares (Table~\ref{tab:fit1}, $5$ degrees of
freedom). The reported $\chi^2$ $p$-values are approximate lack-of-fit
diagnostics using bootstrap-based standard errors, not formal tests of an
underlying asymptotic law. A one-parameter decay $c/\ln x$ is adequate only for the $K_n > 0$ shift. It is
rejected for $\gcd = 2$, $K_n = 0$, and $K_n < 0$, and the alternative
$c(\ln\ln x)/\ln x$ is rejected much more strongly in those three columns, and
is disfavored for $K_n > 0$ ($p = 0.20$). In other words, the $\gcd = 2$, $K_n = 0$, and
$K_n < 0$ shifts decay somewhat faster than $1/\ln x$ over this range, whereas the
$K_n > 0$ shift decays like $c/\ln x$ with $c \approx -4.6$.

\begin{table}[htbp]
\centering
\small
\begin{tabular}{lrrrrrr}
\toprule
 & \multicolumn{3}{c}{$c/\ln x$} & \multicolumn{3}{c}{$c\,(\ln\ln x)/\ln x$} \\
\cmidrule(lr){2-4}\cmidrule(lr){5-7}
Shift & $c$ (pp) & $\chi^2$ & $p$ & $c$ (pp) & $\chi^2$ & $p$ \\
\midrule
$\gcd = 2$ & $+30.44 \pm 0.04$ & 122 & $1\times10^{-24}$ & $+9.72 \pm 0.01$ & 598 & $5\times10^{-127}$ \\
$K_n > 0$ & $-4.58 \pm 0.02$ & 0.55 & 0.99 & $-1.46 \pm 0.01$ & 7.3 & 0.20 \\
$K_n = 0$ & $-3.26 \pm 0.01$ & 281 & $9\times10^{-59}$ & $-1.04 \pm 0.00$ & 540 & $2\times10^{-114}$ \\
$K_n < 0$ & $+7.84 \pm 0.03$ & 27 & $6\times10^{-5}$ & $+2.51 \pm 0.01$ & 107 & $2\times10^{-21}$ \\
\bottomrule
\end{tabular}
\caption{One-parameter weighted least-squares fits to the six interval shifts
($5$ degrees of freedom).}\label{tab:fit1}
\end{table}

Two-parameter forms fit all four shifts (Table~\ref{tab:fit2}, $4$ degrees of freedom):
a power law $c/(\ln x)^a$, and $c/\ln x + d/(\ln x)^2$. The fitted exponents are $a \approx
1.25$ ($\gcd=2$), $0.95$ ($K_n>0$), $1.78$ ($K_n=0$), and $1.31$ ($K_n<0$). These fits
are exploratory: the parameters are strongly correlated and individually poorly
determined, and the two functional forms cannot be distinguished by the data
($p$-values between $0.28$ and $1.00$). We therefore draw only the following conclusions
over $10^5 \le x \le 2\times 10^{10}$: the shifts decrease monotonically in magnitude;
they are compatible with decay to zero; and the data do not determine the asymptotic
rate.

\begin{table}[htbp]
\centering
\small
\begin{tabular}{lrrrrrr}
\toprule
 & \multicolumn{3}{c}{$c/(\ln x)^a$} & \multicolumn{3}{c}{$c/\ln x + d/(\ln x)^2$} \\
\cmidrule(lr){2-4}\cmidrule(lr){5-7}
Shift & $(c, a)$ & $\chi^2$ & $p$ & $(c, d)$ & $\chi^2$ & $p$ \\
\midrule
$\gcd = 2$ & $(66.6,\ 1.25)$ & 4.92 & 0.30 & $(23.4,\ 161.9)$ & 5.04 & 0.28 \\
$K_n > 0$ & $(-3.97,\ 0.95)$ & 0.36 & 0.99 & $(-4.76,\ 4.09)$ & 0.38 & 0.98 \\
$K_n = 0$ & $(-37.9,\ 1.78)$ & 1.69 & 0.79 & $(-0.74,\ -57.5)$ & 2.06 & 0.73 \\
$K_n < 0$ & $(20.9,\ 1.31)$ & 0.19 & 1.00 & $(5.51,\ 53.1)$ & 0.08 & 1.00 \\
\bottomrule
\end{tabular}
\caption{Two-parameter fits to the six interval shifts ($4$ degrees of freedom). Parameters
are point estimates; they are strongly correlated and their uncertainties are
not quoted.}\label{tab:fit2}
\end{table}

\section{Machine verification in Lean~4}\label{sec:lean}

The finite and algebraic results in Sections~\ref{sec:construction}--\ref{sec:pick} were
formalized in Lean~4 with Mathlib, centered on two main modules
(\texttt{PrimeGaps} and \texttt{PrimeTrapezoidPick}) and supported by additional
verification modules. Both modules index primes from $0$ via
\texttt{p := Nat.nth Nat.Prime} (\texttt{p 0 = 2}) and define
\texttt{g n := p (n + 1) - p n}, so the paper's $p_n$ corresponds to \texttt{p (n - 1)};
their definitional equality is verified by \texttt{p\_eq} and \texttt{g\_eq} in
\path{PrimeGaps/TrapezoidAxiomChecks.lean}.

In \texttt{PrimeTrapezoidPick}, the boundary count \texttt{B} is defined by
\eqref{eq:Bn}, while the strict interior lattice points of the translated trapezoid
$Q_n - (p_n, p_n)$ (with vertices $(0,0)$, $(a,a)$, $(a,a+b)$, $(0,a)$ for $a = g_n$,
$b = g_{n+1}$) are defined directly as a finite set:
\[
\texttt{interiorPoints a b} = \bigl\{(u,v) \in \N^2 : 0 < u < a,\; u < v,\;
a v < a^2 + b u\bigr\}.
\]
Note that $0 < u < a$ and $a v < a^2 + b u$ automatically imply $v < a + b$, which
justifies restricting $v$ to $\texttt{range (a + b)}$ in the vertical-slice lemma
\texttt{card\_slice}. The theorem \path{two_mul_card_interiorPoints} proves directly
that $2 \cdot |\texttt{interiorPoints}(g_n, g_{n+1})| = 2I_n$ for all $n \ge 0$ via
vertical slicing (\texttt{card\_slice}) and Gauss pairing (\texttt{floor\_pair}), without
assuming a general library version of Pick's theorem. Thus Proposition~\ref{prop:pick} as
a geometric lattice-point count is verified directly for $I_n$, whereas the geometric
boundary count \eqref{eq:Bn} is taken as a definition (Pick's theorem itself is not
formalized). Concrete kernel checks \path{card_interiorPoints_2_4} and
\path{card_interiorPoints_4_2} (discharged by \texttt{decide}, without
\texttt{native\_decide}) verify $|\texttt{interiorPoints } 2\; 4| = 2$ and
$|\texttt{interiorPoints } 4\; 2| = 7$, matching $Q_3$ and $Q_4$ in
Figure~\ref{fig:trapezoids}.

\begin{table}[htbp]
\centering
\footnotesize
\begin{tabular}{@{}p{0.31\linewidth}p{0.62\linewidth}@{}}
\toprule
Paper statement & Lean~4 declaration(s) \\
\midrule
Lemma~\ref{lem:identity} & \texttt{lemma1}, \texttt{delta\_alt}, \texttt{sum\_delta\_eq} \\
Proposition~\ref{prop:basic} & \texttt{delta\_pos}, \texttt{gap\_even}, \texttt{delta\_even}, \texttt{delta\_ge\_four} \\
Theorem~\ref{thm:sandwich} & \texttt{sum\_delta\_ge\_half}, \texttt{theorem2\_corrected}, \texttt{sum\_delta\_lower} \\
Theorem~\ref{thm:min_I} & \texttt{twoI\_zero}, \texttt{twoI\_one}, \path{twoI_of_g_eq_two}, \\
 & \path{not_both_gaps_two}, \texttt{twoI\_ge\_two}, \texttt{twoI\_ge\_four}, \\
 & \path{twoI_ge_fourteen_of_g_ge_four}, \texttt{twoI\_eq\_four\_iff}, \\
 & \path{twoI_eq_fourteen_iff} \\
Theorem~\ref{thm:boundary} (eq.~\eqref{eq:Bsum} only) & \texttt{sum\_B\_eq}, \texttt{sum\_B\_le} \\
Direct interior-point formula & \texttt{interiorPoints}, \texttt{card\_slice}, \texttt{floor\_pair}, \\
 & \path{two_mul_card_interiorPoints} \\
Concrete checks \& module bridge & \path{card_interiorPoints_2_4}, \path{card_interiorPoints_4_2}, \\
 & \texttt{p\_eq}, \texttt{g\_eq} \\
\bottomrule
\end{tabular}
\caption{Correspondence between paper results and verified Lean~4 declarations.}\label{tab:lean}
\end{table}

The declarations in Table~\ref{tab:lean} were checked using \texttt{\#print axioms}.
The command \path{lake build PrimeGaps.TrapezoidAxiomChecks} prints their
axiom dependencies: \texttt{propext}, \texttt{Classical.choice}, and
\texttt{Quot.sound}, with no \texttt{sorry}. Corollary~\ref{cor:pnt}, the asymptotic limit and $O(1/\ln N)$ rate in
Theorem~\ref{thm:boundary}, and the numerical computations in Section~\ref{sec:numerics}
are external to the formalization.

\paragraph{Data and code availability.} The Lean~4 source files,
\texttt{lake-manifest.json}, and two Python scripts—\texttt{reproduce\_tables.py}
(Tables~\ref{tab:small} and~\ref{tab:density}) and \texttt{prime\_gap\_stats.py}
(the resumable streaming script producing Tables~\ref{tab:gcd}--\ref{tab:fit2}, the
convexity counts, and the fits up to $3\times 10^{10}$)—are available at
\url{https://github.com/Bronxproject/lattice-trapezoids-prime-gaps} and archived on
Zenodo at \url{https://doi.org/10.5281/zenodo.24451767}.

\section{Open problems}\label{sec:open}

We conclude with two questions suggested by the computations in
Section~\ref{sec:numerics}.

\begin{enumerate}
\item Under a reduced-residue heuristic motivated by Hardy--Littlewood-type prime-tuple
  equidistribution, the empirical distribution of
  $\gcd(g_n, g_{n+1})$ is predicted to approach the reduced-residue model $d_{2k}$ as $x \to \infty$. Can an extension of the Lemke~Oliver--Soundararajan
  conjecture~\cite{los} to triples of consecutive primes predict both the single-gap
  deficit and the adjacent-gap shift in Table~\ref{tab:gcd}, which for $\gcd = 2$
  decreases from $+2.90$ to $+1.29$ percentage points between $10^5$ and $10^{10}$
  (Table~\ref{tab:scale}) faster than $1/\ln x$?
\item For the discrete convexity statistic $K_n = \tfrac12(g_n g_{n+2} - g_{n+1}^2)$, the
  concavity excess relative to the local windowed marginal null decreases from
  $+0.840$ to $+0.332$ percentage points over $10^5 \le x \le 2\times 10^{10}$. A
  one-parameter $c/\ln x$ decay is rejected ($p = 6\times 10^{-5}$), whereas two-parameter
  forms fit well (Table~\ref{tab:fit2}). Does this excess tend to zero as $x \to \infty$,
  and at what rate?
\end{enumerate}

\paragraph{Acknowledgments.} The author conceived the geometric construction,
formulated the main questions, directed the mathematical and statistical analysis,
and takes full responsibility for the content of this paper. AI assistants (ChatGPT, Claude,
and Gemini) were used as tools to assist in drafting Lean~4 proof scripts, numerical
computation code, and LaTeX prose. Every formal theorem was independently verified
by the Lean~4 kernel with no \texttt{sorry}, and all numerical tables and fits were
generated directly from the accompanying scripts.

\end{document}